\documentclass[a4paper,11pt]{article}
\usepackage[utf8]{inputenc}
\usepackage[french]{babel}
\usepackage{graphicx}

\usepackage{ifpdf}

\ifpdf
\fi

\begin{document}

\thispagestyle{empty}
{\Large
\noindent Mémoire de Formation Musicale 3\ieme\ cycle \hfill Mai 2015

\vspace{1ex}

\noindent Présenté par Nicolas Trotignon
}
\vspace{9cm}

\begin{center}{\Huge\center
Sur le théorème des trois distances

\vspace{1ex}

 et la construction des gammes}
\end{center}

\vspace{9cm}

{\Large
\begin{center} --- \'Ecole Nationale de Musique de Villeurbanne --- \end{center}
}
\newpage 
\thispagestyle{empty}
\mbox{}
\newpage 
\setcounter{page}{1}
\renewcommand{\og}{``}
\renewcommand{\fg}{''}

\rule{9em}{0ex}\parbox{8cm}{{\it\og Ains n’avez-vous paour, amy, que
  tousiours couché comme ung veau et roulant la vastitude de ces
  choses en la sphéréité de vostre entendement, elles ne cataglyptent
  une façon de microsme en votre personne et ne vous appréhendent
  vous-même ?\fg{}}, Gustave Flaubert, Lettre à Louis Bouilhet du 26
  décembre 1852.}

\vspace{2ex}
\section*{Introduction}


Ce mémoire trouve son origine dans une discussion informelle avec
Guillaume Hanrot, professeur d'informatique à l'\'Ecole Normale
Supérieure de Lyon et musicien amateur.  Il a remarqué qu'un théorème,
démontré dans les année 1950 par la mathématicienne Vera
S\'os et connu sous le nom de \og théorème des trois distances\fg{},
permet d'éclairer certains aspects de la théorie des gammes musicales.
Après quelques recherches infructueuses sur internet, j'ai cru que ce
lien n'avait pas encore été découvert par les musicologues.  Mais une
fois décidé le sujet de ce mémoire, des recherches plus poussées m'ont
conduit à la thèse de Norman Carey, soutenue en 1998 à l'université de
Rochester aux Etats-Unis, et qui porte essentiellement sur les
propriétés formelles des gammes et en particulier, leur lien avec
le théorème sus-mentionné.  Je me trouvais d'une part conforté dans
l'idée que le sujet était digne d'intérêt, et d'autre part un peu
déçu d'avoir été précédé. De ce fait, on ne trouvera ici rien ou
presque qui ne soit déjà présent dans la littérature.

On s'intéresse donc ici à l'une des multiples facettes du lien entre
mathématiques et musique.  Ce lien a fait couler beaucoup d'encre
depuis des siècles, et c'est aujourd'hui encore un domaine actif de
recherche.  De par ma formation (j'ai étudié les mathématiques, et
longuement), je me trouvais entraîné comme par une pente vers ce type
de sujets, mais les pentes étant plus avantageusement remontées que
descendues, j'aurais souhaité travailler sur autre chose. J'ajouterais
que je ne regarde pas sans un certain scepticisme le lien entre
mathématiques et musique.  Pour réel et profond qu'il soit, son étude
 ne me semble pas toujours produire la meilleure des
musiques, ni les meilleures des mathématiques.  En outre,
l'éloignement et la technicité des deux sujets qu'il se propose de
rapprocher limitent le nombre des contradicteurs, ce qui est propice à
diverses impostures, et paradoxalement à une certaine nonchalance.
Bref, sans ce point de départ intrigant, à savoir ce magnifique
théorème des trois distances, et sans l'enthousiasme de Claire Lapalu
(professeur de formation musicale et responsable des mémoires de
3\ieme~cycle à l'ENM), j'aurais sûrement choisi de travailler sur un
autre sujet, plus proche de mes intérêts musicaux, sans doute sur la
chanson ou les musiques vocales.

La première partie de ce mémoire explique le problème général de la
construction de gammes.  La deuxième partie donne quelques rappels
indispensables sur les logarithmes.  Dans la troisième partie, le
théorème des trois distances est expliqué et appliqué à la
construction de gammes.  Nous montrons comment il permet de retrouver
naturellement les gammes pythagoricienne, pentatonique, diatonique et
chromatique, ainsi que les demi-tons diatoniques et chromatiques, puis
nous expliquons un procédé général de construction de gammes provenant
de la thèse de Norman Carey.

\subsection*{Remarques préalables}

Avant de poursuivre, précisons que ce travail est conçu pour être lu
par des musiciens ayant un bagage mathématique minimal (disons lycée
ou fin de collège).  Les notions de solfège comme \og intervalle\fg{},
\og quinte\fg{}, \og transposition\fg{}, etc sont donc supposées
connues du lecteur, au contraire des notions mathématiques, notamment
les logarithmes, qui sont expliquées en détail.  Notons que les
logarithmes sont inutiles pour énoncer, comprendre, ou même prouver le
théorème des trois distances.  C'est son utilisation en musique qui
nécessite les logarithmes.

Dans ce travail, on distingue la construction des gammes du problème
bien connu des musiciens sous le nom de \og tempérament\fg{}.
Expliquons ce point.  Une fois une gamme plus ou moins fixée (par
exemple do-ré-mi-fa-sol-la-si-do), le tempérament consiste en
l'ajustement exact des hauteurs des notes en fonction de différentes
contraintes (consonnance des intervalles, possibilité de transposer
une mélodie ou de moduler, fidélité à la musique ancienne, recherche
d'une couleur, facture et accord d'instruments, etc).  Le problème du
tempérament est selon moi bien plus musical que mathématique.

Le problème de la construction de gammes est plus général: on ne
présuppose aucune échelle, même approximative.  En particulier le
nombre de notes de la gamme n'est pas connu d'avance: il n'y a pas
{\it a priori} 5 notes (gamme pentatonique), 7 notes (gamme
diatonique) ou 12 notes (gamme chromatique).  Le problème de la
construction de gammes ne se pose pas de manière urgente aux musiciens
qui ne semblent pas réclamer instamment de nouvelles gammes (alors
qu'ils réclament que leurs instruments soient accordés, ce qui
implique le choix d'un tempérament).  Ainsi n'y a-t-il pas
d'inconvénient majeur à considérer la question sous un angle
exclusivement mathématique.  Le problème est au demeurant surtout
retrospectif, en ce qu'il éclaire le pourquoi des gammes aujourd'hui
utilisées.  Dans ses aspects exploratoires, il relève de la pure
spéculation, avec tout ce que cela implique de promesse vague et de
risque.

On se donne donc certaines contraintes musicales exprimées
mathématiquement, et on cherche (mathématiquement encore) l'ensemble
des sons permettant au mieux de les satisfaire.  Une des conclusions de
ce travail est que les nombres 5, 7, et 12 qu'on rencontre dans de
nombreuses gammes ne proviennent pas seulement de contingences
culturelles, mais ont des propriétés mathématiques remarquables qui
expliquent {\it a posteriori} leur présence.  Ce fait est connu depuis
longtemps, mais le théorème des trois distances l'illustre de manière
claire et saisissante.  Cependant, le théorème des trois distances
n'éclaire pas le problème du tempérament tel que défini ci-dessus.

\subsection*{Sur la bibliographie}

Sur toutes les généralités concernant la construction des gammes,
notamment le procédé pythagoricien et la définition des différents
commas, Wikipédia est une source remarquablement claire (et fiable
pour autant que je puisse en juger). La plupart des informations
générales de ce mémoire en proviennent.  Sur le tempérament, le livre de
Pierre-Yves Asselin~\cite{asselin:mt} m'a paru tout à fait
remarquable.  Il est principalement centré sur la résolution pratique
de véritables problèmes musicaux, comme l'accordage des clavecins.  Sa
lecture est donc un peu frustrante dans la mesure où elle ne donne pas
à \emph{entendre} l'effet musical des différents tempéraments évoqués.
Dans le livre de Gareth Loy~\cite{Loy:musimathicsV1}, on trouve une
bonne synthèse sur les mathématiques et la musique, qui explique bien
les bases de l'accoustique et de la construction des gammes.  Une
preuve du théorème des trois distances peut être consultée dans le
livre de Jean-Paul Allouche et Jeffrey
Shallit~\cite{DBLP:books/daglib/0025558}. On peut aussi consulter la
synthèse de Pascal Alessandri et Valérie
Berthé~\cite{alesandriniBe:98}.  Sur le lien entre ce théorème et la
musique, les seules sources que j'ai trouvées sont la thèse de Norman
Carey~\cite{carey:these} et l'article qu'il a co-signé avec David
Clampitt~\cite{careyCl:2012}.

\section{Construction de gammes}

Un \emph{son} est la perception que l'oreille humaine a d'une
oscillation périodique de la pression de l'air, appelée \emph{onde
  sonore}.  La \emph{hauteur} d'un son s'exprime mathématiquement
par sa fréquence, exprimées en hertz, c'est-à-dire par le  nombre
d'oscillations par seconde de l'onde qui le produit. Pour fixer les
idées concrètement, un son de $440$ hertz est un La obtenu juste à
droite de la serrure d'un piano.  L'oreille humaine humaine peut
entendre des sons d'environ 20\,Hz (le plus grave) à 20\,000\,Hz (le plus
aigu).  Entre ces bornes extrêmes, l'ensemble des hauteurs possibles
est un continuum, c'est-à-dire que n'importe quel nombre
réel\footnote{Un nombre réel est un nombre \og avec ce qu'on veut
  après la virgule\fg{}.} est la fréquence d'un son.  Il y a donc
potentiellement une infinité de hauteurs possibles.

Pour des raisons culturelles et pratiques (touches d'un clavier,
frettes d'un manche, trous d'une flûte\dots), les musiciens préfèrent
limiter les hauteurs de sons à un ensemble fini. Construire une gamme,
c'est choisir les fréquences des sons auxquels on se restreint.  Le
problème peut paraître abstrait, mais il s'est posé concrètement, bien
qu'en des termes différents, aux tout premiers facteurs et accordeurs
d'instruments, sans doute à des facteurs de flûtes de la préhistoire
(la page Wikipédia \og instrument de musique\fg{} mentionne les restes
d'une flûte en os fabriquée 45\,000 ans avant notre ère).  La théorie
des ondes sonores date seulement du XVII\ieme\ siècle avec les travaux
de Mersenne et Galilée, mais il se trouve que la fréquence d'un son
émis par une corde grattée est inversement proportionnelle à la
longueur de la corde, ce qui a permis sans théorie ni instrument de
mesure sophistiqué de quantifier précisemment la hauteur des sons dès
la plus haute antiquité.  Les aspects mathématiques du problème de la
construction des gammes sont donc reconnus depuis au moins 2\,500 ans,
et les Grecs anciens attribuent cette découverte à l'école de
Pythagore.  On en vient donc au premier principe de la construction de
gammes (que l'on choisit d'exprimer mathématiquement).

\vspace{1ex}
\noindent{\bf Principe 1:} une gamme est un ensemble fini (et en
pratique assez petit) de nombres, qui sont les fréquences des
différents sons de la gamme.

\vspace{1ex}

Les pythagoriciens, ou peut-être Pythagore lui-même (on manque de
sources sur ces temps reculés), ont remarqué qu'en doublant la
fréquence d'un son, on obtient un deuxième son très consonnant avec le
premier, à tel point que, dans certains contextes musicaux, ces deux
sons peuvent être considérés comme équivalents (dans notre système
musical, il s'agit d'une note et de la même note à l'octave
supérieure).  Ceci mène au deuxième principe de la construction des
gammes, le principe de l'identité des octaves, explicité pour la
première fois par Jean-Philippe Rameau dans son {\it Traité de
  l'harmonie réduite à ses principes naturels}, mais bien évidemment
utilisé implicitement par les musiciens l'ayant précédé:

\vspace{1ex}

\noindent{\bf Principe 2:} si un nombre est dans la gamme, le nombre
obtenu en le doublant est considéré comme la fréquence d'un son
équivalent musicalement. 

\vspace{1ex}

De même, en mulpliant la fréquence d'un son par d'autres nombres
rationnels simples\footnote{Un nombre rationnel est le résultat de la
  division de deux nombres entiers (\og entier\fg{} signifie sans rien
  après la virgule).}, comme par exemple $3/2$ ou $4/3$, $5/4$ ou
$6/5$, on obtient un nouveau son formant un intervalle assez
consonnant avec le premier, et perçu par l'oreille humaine comme
musicalement intéressant.  Par exemple, en multipliant par $3/2$, on
obtient ce qu'on appelle dans notre système musical une quinte pure,
par $4/3$ une quarte pure, par $5/4$ une tierce majeure pure et par
$6/5$ une tierce mineure pure.  Notons toutefois qu'à ce stade, les
mots \og tierce\fg{}, \og quarte\fg{} et \og quinte\fg{} sont
fallacieux puisqu'ils présupposent un écart de $3$, $4$ ou $5$ notes
entre les deux sons, alors qu'aucune gamme (et aucun nombre de notes)
n'est présupposé.  Il s'agit simplement d'intervalles consonnants,
qu'il faudrait peut-être nommer autrement.

Les auteurs se perdent en conjectures sur l'origine de cette
perception de consonnance, à tel point qu'on peut se demander si la
conclusion de tel ou tel n'en dit pas plus sur lui-même que sur
l'origine réelle de la consonnance.  Est-elle physiologique, acoustique,
culturelle, numérologique, mystique ou mathématique?  Y a-t-il
d'ailleurs une origine \og réelle\fg{} de la consonnance?  Pour la
commodité de l'écriture de ce mémoire, nous tenons la consonnance des
intervalles provenant de rapports rationnels simples pour un fait
brut que nous ne cherchons pas à expliquer plus avant.  On en vient au
troisième principe de la construction des gammes:

\vspace{1ex}

\noindent{\bf Principe 3:} si un nombre est dans la gamme, on souhaite
que les nombres obtenus par multiplication par des rapports
simples, comme $2$, $3/2$ ou $4/3$ (correspondant à des intervalles
consonnants ascendants), et leurs inverses $1/2$, $2/3$ ou $3/4$
(correspondant aux mêmes intervalles descendants), soient aussi présents
dans la gamme.

\vspace{1ex}

Pris ensemble, les trois principes sont mathématiquement
contradictoires. On peut prouver qu'il est impossible de construire un
ensemble fini de nombres réels clots par multiplication par $1/2$, $2$,
$3/2$ et $2/3$.  La preuve n'est pas très compliquée, mais son intérêt
pour ce mémoire étant limité nous l'omettons.  Indiquons qu'elle se
ramène entre autres au fait qu'en multipliant indéfiniment des $3$, on
n'obtiendra jamais un multiple de $2$.  L'impossibilité de construire
une gamme pleinement satisfaisante vient donc d'une impossibilité
mathématique fondamentale, et non d'un manque d'habileté des
mathématiciens ou des musiciens.  C'est de là que la construction de
gammes est réellement un \emph{problème}.

Devant l'impossibilité de construire une gamme répondant pleinement
aux trois principes, plusieurs stratégies sont envisageables: on peut
s'affranchir partiellement du premier principe en acceptant qu'il y
ait un nombre de sons possibles, non pas infini, mais s'en approchant,
c'est-à-dire très grand (cette méthode est par exemple utilisée dans
le tempérament à 53 intervalles égaux de Holder qu'on rencontrera dans
la troisième partie de ce mémoire); on peut se contenter
d'une approximation des rapports simples du deuxième principe (c'est
la solution actuellement retenue dans ce qu'on appelle le tempérament
égal); on peut garder des rapports exacts pour un maximum
d'intervalles, en en \og sacrifiant\fg{} certains autres qui
deviennent inutilisables musicalement (c'est la solution retenue dans
de nombreux systèmes anciens de tempérament, où une quinte est \og
sacrifiée\fg{}, la fameuse quinte du loup).

\section{Logarithmes}

Il est hélas difficile de poursuivre sans un peu de bagage
mathématique.  L'une des difficultés techniques pour comprendre le
lien entre le théorème des trois distances et la construction des
gammes est que l'accoustique suit une logique multiplicative (on
obtient un intervalle de quinte pure en \emph{multipliant} par
$3/2 = 1.5$), tandis que les musiciens pensent plutôt les intervalles
de manière additive, ce en quoi ils ont parfaitement raison (on
obtient une quinte en \emph{additionnant} une tierce majeure et une
tierce mineure, sûrement pas en les multipliant).  Ce passage d'une
logique multiplicative à une logique additive est tout à fait
légitime. Les mathématiciens l'appellent un \og passage au
logarithme\fg{}, et je n'ai malheureusement pas trouvé de moyen de me
passer de cette notion (le mathématicien Henri Cartan disait : \og Il
est plus facile d’apprendre les mathématiques que d’apprendre à s’en
passer\fg{}).

Le logarithme est une notion tout à fait concrète, disons aussi
concrète que l'addition ou la multiplication.  Tout d'abord, le
logarithme est une \emph{fonction}, c'est-à-dire que c'est un procédé
qui transforme un nombre en un autre nombre.  Par exemple, le
logarithme de 1 vaut 0, ce qui se note $\log(1) = 0$.  Le logarithme
de $12.7$ est un nombre assez difficile à calculer à la main; on ne
peut même pas donner de niveau scolaire correspondant à cette tâche,
puisque le calcul numérique a presque entièrement disparu des
programmes de l'enseignement secondaire.  Mais la plupart des
calculatrices, utilisées à partir du collège et responsables de cette
disparition, ont pour se la faire pardonner une touche log qui permet
de trouver que le logarithme de $12.7$ vaut environ $1.1038$, ce qui
s'écrit $\log(12.7) = 1.1038$.

Voilà donc un problème de résolu: grâce à une calculatrice, on peut
tout à fait concrètement calculer le logarithme de n'importe quel
nombre.  Mais que signifie le logarithme? Pour s'en faire une première
idée approximative, il suffit de compter les chiffres du nombre, puis
d'enlever $1$.  Par exemple, $103$ s'écrit avec $3$ chiffres (un, zéro
et trois), son logarithme est donc à peu près $3-1=2$, en fait
$2.0128$, soit $2$ \og et des poussières\fg{} comme on peut le
vérifier à la calculatrice.  Si on s'approche de nombres à $4$
chiffres (donc de $1000$), l'approximation devient de moins moins
bonne, par exemple $\log(998) = 2.9991$, ce qui est assez \og
éloigné\fg{} de $2$, et presque égal à $3$ (998 a \og presque\fg{} 4
chiffres).  Le nombre $367\,987\,434$ a $9$ chiffres, son logarithme sera
environ de $9-1=8$ (et encore des poussières).  Cette définition
approximative (car on ne dit pas comment calculer les \og
poussières\fg{}, qui seront pourtant essentielles par la suite)
permet de saisir quelques propriétés importantes du logarithme.  Par
exemple, qu'un nombre très grand a un logarithme assez petit: un
nombre dont le logarithme serait $50$ doit être vraiment gigantesque,
puisqu'il s'écrit avec $51$ chiffres.  Et aussi, que quand on
multiplie un nombre par $10$, on lui ajoute un chiffre, ce qui fait
que son logarithme augmente de $1$.

Pour les nombres $1$, $10$, $100$, etc, (ce qu'on appelle les
puissances de $10$) la méthode de calcul du logarithme donnée ci-dessus
est en fait exacte, et le logarithme est précisement le nombre de
chiffres moins~$1$, ou ce qui revient au même, le nombre de zéros.  Et
il est bien connu que pour multiplier des puissances de $10$, il
suffit d'additionner le nombre de zéros (par exemple, $100$ a $2$
zéros, $1000$ a $3$ zéros, et $100\times 1000 = 100000$ a $2+3 = 5$
zéros).  Ainsi, quand on mulptiplie des nombres, leurs logarithmes
s'additionnent (du moins pour les puissances de $10$).  Cela s'écrit
$\log(a \times b) = \log(a) + \log(b)$.  Cette formule ci-dessus est
en fait vraie pour tous nombres réels strictement positifs $a$ et $b$.
La définition approximative du logarithme donnée ici ne permet pas de
prouver cette relation.  En effet, calculer ou définir précisement le
logarithme est difficile (il faut ajuster ces fameuses \og
poussières\fg{} en s'arrangeant pour que la formule soit vraie), et
l'existence même des logarithmes est un fait mathématique tout à fait
remarquable découvert seulement au XVII\ieme\ siècle, c'est-à-dire
bien après la théorie du tempérament qui date de l'Antiquité.  Les
logarithmes peuvent être vus comme un dictionnaire qui traduit les
multiplications en additions.  Ce fut d'ailleurs jusqu'à l'avènement
des calculatrices l'une de leurs applications: pour mutiplier
deux nombres (opération fastidieuse), on calculait leur logarithme
grâce à une table établie une fois pour toute, on additionnait les
résultats (opération facile), puis on regardait à nouveau dans la
table pour voir de quel nombre le résultat était le logarithme.  Ce
nombre est le résultat de la multiplication.  L'addition pouvait
d'ailleurs être effectuée analogiquement, c'est-à-dire en mettant
bout-à-bout deux réglettes graduées de longueur égale aux nombres à
additionner, et en trouvant le résultat par une mesure de la longueur
totale: c'est le principe de la fameuse règle à calcul, aujourd'hui
tombée en désuétude.

Dans le cadre de ce travail, nous admettons simplement qu'il existe
une fonction notée $\log$, dont l'existence est matérialisée par la
touche \og log\fg{} des calculatrices, et telle que pour tous nombres
$a$ et $b$, on a: $\log(a\times b) = \log(a) + \log(b)$.

Pour ne pas compliquer la discussion, j'ai omis jusque-là de préciser
qu'il y a plusieurs fonctions logarithmes différentes.  Celle que j'ai
approximativement définie ici est ce que les mathématiciens appellent
le logarithme en base $10$. Il a la propriété que lorsqu'on multiplie un
nombre par $10$, son logarithme augmente de $1$.  Pour la théorie des
gammes, la multiplication par $10$ n'a pas grand intérêt, alors que
la multiplication par $2$ est riche de sens (elle fournit l'intervalle
d'octave).  On utilisera donc plutôt le logarithme en base $2$.  Il peut
être obtenu avec une calculatrice en calculant le logarithme en base
$10$, puis en le multipliant par une constante à peu près égale à $3.32$.  Une
propriété intéressante du logarithme en base $2$ est que si un nombre
est multiplié par $2$, alors son logarithme augmente de $1$.

\`A partir de maintenant, $\log(x)$ désigne le logarithme en base $2$ de
$x$.  Il est très important de remarquer que $\log(2) = 1$.  Et la
formule $\log(a\times b) = \log(a) + \log(b)$ reste vraie.

Nous sommes maintenant mûrs pour reformuler le troisième principe de
la construction des gammes.  Comme la multiplication est une opération
plus difficile que l'addition, il est plus simple de rechercher non
pas les fréquences des sons, mais leur logarithme.  On choisit
certains intervalles simples, comme la quinte pure, ou la tierce pure,
sous la forme de rapports rationnels, $3/2$ et $5/4$ comme déjà
expliqué.  On souhaite que, si une note est dans la gamme, alors en
additionnant (et non plus en multipliant) un nombre correspondant à
l'octave ou à la quinte pure par exemple, on reste dans la gamme.  Ces
nombres à additionner ne sont autres que les logarithmes des rapports
précédents.  L'élévation d'une octave correspond donc à l'ajout de 1
(car $\log(2)=1$).  On le verra dans la troisième partie, le nombre
$0.58496$ joue un grand rôle dans la construction de la gamme
pythagoricienne.  C'est qu'il n'est autre que le logarithme de $3/2$.

Ce qu'on peut retenir de ce chapitre, c'est qu'ajouter $0.58496$ à la
fréquence d'un son (ou plutôt à son logarithme) est une opération
naturelle musicalement (il s'agit de s'élever d'une quinte).  Il est
important de comprendre que ce nombre n'a rien d'arbitraire, du moins
n'est-il pas plus arbitraire que la quinte elle-même. Il ne provient
ni des mathématiques (même si les mathématiques participent à une
étape purement technique de son élaboration), ni d'une quelconque
mystique numérologique, comme si l'on avait choisi $\pi$, ou la
racine carrée de l'âge du chef d'orchestre, ou le nombre d'or
(quoiqu'étudier le nombre d'or dans ce contexte ne serait peut-être
pas si absurde étant données les propriétés remarquables de son
développement en fraction continue, mais cela sort du cadre de ce
travail).

\section{Le théorème des trois distances}

Le théorème des trois distances concerne l'addition répétitive d'un
même nombre, ce qui est intéressant musicalement, puisque comme on l'a
vu, cela peut être interprété comme l'ajout successif de quintes
ou d'autre intervalles.  Il ne provient pas de recherches sur la
musique et les sources à propos de son application en musicologie ne
remontent pas au-delà de la fin des années 1990.  J'estime toutefois
probable que certains mathématiciens bien au fait du procédé pythagoricien
de constructions des gammes ont pu (à l'instar de Guillaume Hanrot)
deviner le lien.  Il me semble plus improbable que des musicologues
professionnels aient deviné ce lien, car si le procédé pythagoricien
est assez largement connu (y compris de nombreux mathématiciens qui
pratiquent la musique), le théorème des trois distances
n'est quant à lui connu que de mathématiciens spécialisés.  Mais
énonçons le théorème.

Postulons un cercle de longueur 1 avec un point de départ. Par
exemple, la piste circulaire d'un stade longue de 1 kilomètre, mais
peu importe l'unité de longueur utilisée:


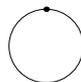
\begin{center}
\includegraphics[width=1cm]{genLog52Pointe.0}
\end{center}


Un marcheur avance à partir du point de départ d'une longueur fixée,
qu'on appelle la \og longueur de son pas\fg{}, par exemple $0.58496$
(nombre intéressant on l'a vu). Le marcheur arrive à un certain
endroit:

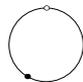
\vspace{-2ex}

\begin{center}
\includegraphics[width=1cm]{genLog52Pointe.1}
\end{center}

\vspace{-3ex}

Il refait un pas de même longueur, puis encore un, et ainsi de suite,
un grand nombre de fois, appelé le \og nombre de pas\fg{}.  On garde
trace sur le cercle de toutes ses stations, mais on ne tient aucun
compte de l'ordre dans lequel les stations ont été visitées.  On se
retrouve donc avec un cercle et un certain nombre de points sur ce
cercle.  Ci-dessous, les cercles successifs après douze pas:

\noindent
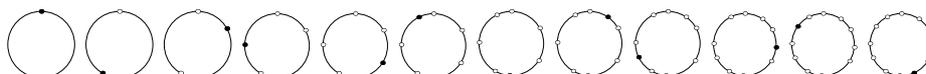
\begin{center}
\includegraphics[width=.9cm]{genLog52Pointe.0}
\includegraphics[width=.9cm]{genLog52Pointe.1}
\includegraphics[width=.9cm]{genLog52Pointe.2}
\includegraphics[width=.9cm]{genLog52Pointe.3}
\includegraphics[width=.9cm]{genLog52Pointe.4}
\includegraphics[width=.9cm]{genLog52Pointe.5}
\includegraphics[width=.9cm]{genLog52Pointe.6}
\includegraphics[width=.9cm]{genLog52Pointe.7}
\includegraphics[width=.9cm]{genLog52Pointe.8}
\includegraphics[width=.9cm]{genLog52Pointe.9}
\includegraphics[width=.9cm]{genLog52Pointe.10}
\includegraphics[width=.9cm]{genLog52Pointe.11}
\end{center}


On regarde maintenant un deuxième marcheur, qui souhaiterait parcourir
les différentes stations du premier, mais cette fois dans leur ordre
naturel sur le cercle. Il devra faire des pas assez petits (et ce
d'autant plus que le premier marcheur a créé beaucoup de stations, qui
\og encombrent\fg{} le cercle), et pas nécessairement tous de la même
longueur sur le cercle. Le théorème des trois distances affirme:

\vspace{2ex}

\noindent{\bf Théorème:} 
Pour toute longueur de pas choisie au départ, et pour tout nombre de
pas du premier marcheur, le deuxième marcheur  rencontrera au
maximum  trois longueurs de pas différentes lors de son parcours.

\vspace{2ex}

Un lecteur méticuleux pourra vérifier la validité de l'énoncé sur
chacun des douze cercles de la figure ci-dessus.  Le théorème concerne
la distance curviligne (c'est-à-dire le long du cercle) entre deux
points, mais celle-ci étant en correspondance bi-univoque avec la
longueur de la corde, un double-décimètre permet de procéder
pratiquement aux vérifications.  Parfois, il y a bien trois distances
différentes, mais parfois seulement deux, nous y reviendrons (le
théorème dit bien trois \emph{au maximum}).  Bien qu'assez
simple, cet énoncé n'est pas évident à démontrer rigoureusement.  Il a
d'abord été conjecturé par Hugo Steinhaus, puis démontré par Vera
S\'os dans les années 1950, avant que plusieurs mathématiciens n'en
donnent diverses preuves.  Il est appelé de différentes manières, ce
qui ne facilite pas les premières recherches à son sujet: théorème de
Steinhaus, conjecture de Steinhaus, théorème des trois longueurs, ou
\og three gap theorem\fg{} en anglais.  \`A son sujet, on peut
consulter le livre de Allouche et Shalit, dont la lecture nécessite de
bonnes connaissances mathématiques (disons bac+2 ou 3).

Le cercle de longueur $1$ représente toutes les fréquences des sons
possibles (ou plutôt, leur logarithme).  Le point de départ représente
une note choisie arbitrairement comme début de la gamme.  Ce peut être
un Do, ou le La 440.  Ma préférence pour établir un lien avec le
solfège classique va au Fa, parce que c'est le point de départ du
cycle des quintes dans le fameux ordre des dièses.  Sur les figures 
suivantes, je pars donc du Fa.  La forme circulaire est pertinente en
raison du principe d'identité des octaves: quand on monte d'une
octave, c'est-à-dire quand on multiplie par $2$ une fréquence, ou qu'en
logarithme on ajoute $1$ (puisque $\log(2) = 1$), on revient à son point
de départ.  Le cercle a cette même propriété: si on avance de $1$, qui
est justement la circonférence du cercle, on revient à son point de
départ.

Monter d'une quinte pure, c'est faire un pas de longueur
$\log(3/2) = 0.58496$.  Si le premier marcheur est un pythagoricien,
il fait donc douze pas, qui correspondent à douze élévations
consécutives d'une quinte, parfois réduite à l'octave (chaque fois
qu'on repasse par le départ).  Le second marcheur lui parcourra les
douze sons dans l'ordre des aiguilles d'une montre, qui forment ce
qu'on appelle la gamme chromatique pythagoricienne.  Un analyse fine
des longueurs de pas du second marcheur montre qu'il parcourt des pas
correspondant à des demi-tons, mais de deux longueurs différentes.
Comme on le voit sur la figure suivante, chaque ton se divise
 en deux demi-ton inégaux, et les demi-tons diatoniques sont bien
inférieurs aux demi-tons chromatiques:

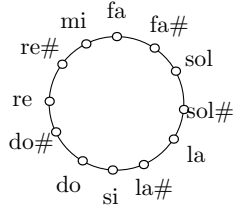
\vspace{-2ex}

\begin{center}
\includegraphics[width=3cm]{cycleQ.11}
\end{center}

\vspace{-3ex}

Le musicien trouvera là une origine à la distinction entre demi-ton
diatonique et demi-ton chromatique. Le mathématicien quant à lui
remarquera que le théorème des trois distances est bien pessimiste: il
prévoit trois distances au maximum pour le second marcheur, or nous n'en
observons que deux, ce qui est souhaitable (deux sortes de demi-ton,
c'est déjà compliqué, trois ce serait bien pire).  En fait cette
situation désirable se rencontre un certain nombre de fois lors du
parcours du premier marcheur.  On la rencontre après $2$ pas, après $5$
pas, après $7$ pas, puis enfin après $12$ pas, comme le lecteur peut le
vérifier sur les quatre figures suivantes :

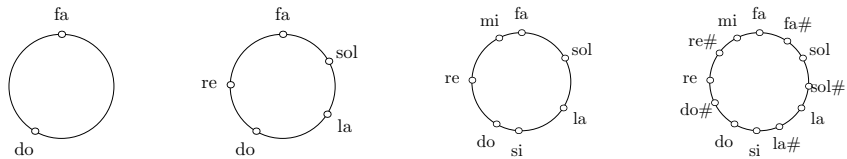
\vspace{-2ex}

\begin{center}
\includegraphics[height=2cm]{cycleQ.1}\hspace{1cm}
\includegraphics[height=2cm]{cycleQ.4}\hspace{1cm}
\includegraphics[height=2cm]{cycleQ.6}\hspace{1cm}
\includegraphics[height=2cm]{cycleQ.11}
\end{center}


Voilà une explication remarquable à la prédominance des gammes à $5$,
$7$ ou $12$ sons dans plusieurs systèmes musicaux.  Cette explication
n'a rien de neuf au demeurant, elle est mathématiquement équivalente à
des considérations bien connues (quoique assez techniques) sur le
développement en fractions continues de $\log(3/2)$, mais le cercle et
le théorème des trois distances en donnent une belle illustration.  Le
procédé est général: on peut partir d'un intervalle quelconque jugé
intéressant, on calcule son logarithme qui donne une longueur de pas,
puis on fait marcher le premier marcheur un certain nombre de fois, en
espérant rencontrer une situation favorable où le second marcheur n'a
que deux distances à parcourir.  Ce procédé fait jouer un rôle central à
l'octave (car on utilise le logarithme en base $2$), mais on peut par
exemple faire des pas égaux à une tierce majeure, et calculer un
logarithme en \og base de quinte\fg{}, c'est-à-dire en base $3/2$ (ce
qui revient à considérer que deux notes séparées d'une quinte sont \og
équivalentes\fg{}), ou dans la base que l'on souhaite.  Norman Carey a
étudié dans sa thèse ce principe général de construction de gammes;
il a montré qu'il conduisait toujours à ce qu'il appelle des gammes \og
bien formées\fg{} (well-formed scales), un concept qu'il a développé à
l'origine indépendamment du théorème des trois distances.  Une
propriété centrale d'une gamme bien formée est que l'intervalle qui
sert à la générer (ici la quinte), se divise toujours en un nombre
égal de \og petites distances\fg{} du second marcheur, ici sept
demi-tons.  Cette propriété, évidente pour le musicien habitué à la
gamme chromatique, est vraie pour toute gamme obtenue selon le procédé
décrit ci-dessus, ce qui n'est pas évident {\it a priori}.  Ceci est
une propriété supplémentaire du théorème des trois distances,
démontrée par Norman Carey (où l'on voit que des questions posées par
des musiciens peuvent faire avancer les mathématiques).

Que se passe-t-il si on ajoute un 13\ieme\ pas ? On tombe sur un son très
proche du Fa (il est bien connu que 12 quintes pures correspondent à
peu près à 7 octaves): 

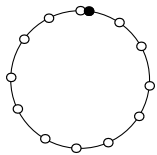
\vspace{-2ex}

\begin{center}
\includegraphics[height=2cm]{genLog52.12}\hspace{1cm}
\end{center}


L'écart entre le Fa et ce 13\ieme\ son est ce qu'on appelle le comma
pythagoricien.  Mathématiquement, après ce 13\ieme\ pas, le deuxième
marcheur se retrouve de nouveau avec trois distances sur le cercle
(dont une minuscule).  En fait, si on retourne à la figure du cercle
de la gamme chromatique avec les noms de notes, on voit que la quinte
qui va du La\# au Fa (ou plutôt de l'enharmonique Sib de La\#,
jusqu'au Fa, pour bien avoir une quinte au sens du solfège, et non une
sixte diminuée), cette quinte donc est inférieure à toutes les autres
quintes du cercle: c'est la seule à contenir $5$ demi-tons diatoniques
au lieu de $4$, et aussi la seule à ne pas être pure: elle ne provient
pas d'une multiplication par $3/2$. C'est la fameuse quinte du loup
des tempéraments anciens.

Mais que se passe-t-il si le premier marcheur marche suffisamment
longtemps, au-delà de ce treizième pas? Se peut-il qu'on se retrouve
de nouveau dans le cas favorable de seulement deux distances?  La
réponse est oui, c'est une autre conclusion du théorème des trois
distances: il est toujours possible de faire marcher plus longtemps le
premier marcheur, jusqu'à un point où le nombre de distances
différentes du deuxième marcheur est au maximum 2.  Le procédé général
fonctionne donc à tous les coups, ce qui n'avait rien de certain
jusque-là, si l'on accepte toutefois de marcher assez longtemps.  Pour
retrouver une telle situation favorable au-delà de 12 pas (dans le cas
de l'élévation par quintes), il faut marcher au total 53 pas. Ceci
mène à la gamme à 53 intervalles de Holder, dont le principe est en
fait connu depuis l'Antiquité (selon Wikipédia, il a été découvert par
le théoricien chinois Ching Fang (78-37 avant JC)).  On conçoit
aisément qu'une gamme à 53 micro-intervalles (égaux au comma de
Holder, qui divise le demi-ton chromatique en 5, et le demi-ton
diatonique en 4) ne soit pas très pratique et que la gamme de Holder
reste cantonnée à la théorie.

On voit donc que la gamme à douze sons est tout à fait remarquable.
Elle a de nombreuses autres propriétés qui ne découlent pas du théorème
des trois distances pour autant qu'on sache.  Par exemple le rapport
entre do et mi (une fois débarrassé des logarithmes, et donc vu à
nouveau multiplicativement) est remarquablement proche de 5/4
(l'intervalle de tierce majeure pure) qui est \og naturellement\fg{},
ou au moins culturellement, consonnant.  Le rapport de ces deux
fréquences est appelé le comma syntonique.  Il vaut $81/80$, ce qui
est remarquablement proche de 1, et je n'ai trouvé nulle part
d'explication à ce phénomène, qui semble une heureuse coïncidence
mathématique.  Autre coïncidence, le rapport entre la tierce mineure
pythagoricienne et 6/5 (tierce mineure pure) est par
chance lui aussi égal au comma syntonique (sans quoi il y aurait un
comma syntonique mineur et un comma syntonique majeur).

Comme on l'a déjà dit, rien n'empêche d'utiliser ce même principe de
construction des gammes en choisissant d'autres valeurs comme point de
départ.  Par exemple, on peut choisir de baser toute la construction sur
l'intervalle de tierce majeure pure 5/4 et les octaves.  Sans préjuger
du résultat, cela peut être motivé par l'harmonie dans la musique
occidentale.  Quels sont les nombres de pas conduisant à des
situations favorables ?  Voici les 29 premières étapes de la
construction:

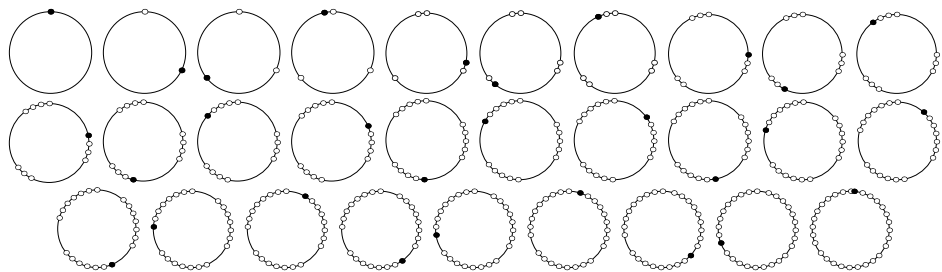
\vspace{1ex}

\begin{center}
\noindent
\includegraphics[width=1.11cm]{genTiercePure.0}
\includegraphics[width=1.11cm]{genTiercePure.1}
\includegraphics[width=1.11cm]{genTiercePure.2}
\includegraphics[width=1.11cm]{genTiercePure.3}
\includegraphics[width=1.11cm]{genTiercePure.4}
\includegraphics[width=1.11cm]{genTiercePure.5}
\includegraphics[width=1.11cm]{genTiercePure.6}
\includegraphics[width=1.11cm]{genTiercePure.7}
\includegraphics[width=1.11cm]{genTiercePure.8}
\includegraphics[width=1.11cm]{genTiercePure.9}
\includegraphics[width=1.11cm]{genTiercePure.10}
\includegraphics[width=1.11cm]{genTiercePure.11}
\includegraphics[width=1.11cm]{genTiercePure.12}
\includegraphics[width=1.11cm]{genTiercePure.13}
\includegraphics[width=1.11cm]{genTiercePure.14}
\includegraphics[width=1.11cm]{genTiercePure.15}
\includegraphics[width=1.11cm]{genTiercePure.16}
\includegraphics[width=1.11cm]{genTiercePure.17}
\includegraphics[width=1.11cm]{genTiercePure.18}
\includegraphics[width=1.11cm]{genTiercePure.19}
\includegraphics[width=1.11cm]{genTiercePure.20}
\includegraphics[width=1.11cm]{genTiercePure.21}
\includegraphics[width=1.11cm]{genTiercePure.22}
\includegraphics[width=1.11cm]{genTiercePure.23}
\includegraphics[width=1.11cm]{genTiercePure.24}
\includegraphics[width=1.11cm]{genTiercePure.25}
\includegraphics[width=1.11cm]{genTiercePure.26}
\includegraphics[width=1.11cm]{genTiercePure.27}
\includegraphics[width=1.11cm]{genTiercePure.28}
\end{center}

\vspace{1ex}

On constate que la situation favorable des deux distances se rencontre
à la troisième étape (ça n'est pas une surprise pour le musicien:
trois tierces majeures valent une octave dans le tempérament égal), et
ensuite à la 28\ieme\ étape seulement.  On obtient une gamme de $28$
micro-intervalles, $3$ longs et 25 petits comme on le voit sur la
figure, presbytie à part.  On pourrait calculer différents commas,
évaluer la fausseté de la dernière tierce, une sorte de \og tierce du
loup\fg{}, mesurer à quel point la gamme permet de s'approcher des
quintes, etc, mais je stoppe là faute de place et de temps.

\section*{Conclusion}

Le théorème des trois distances explique bien certains aspects des
constructions classiques de gammes.  Les figures circulaires de ce
mémoire sont selon moi plus faciles à comprendre que les
considérations classiques sur le développement en fraction continue du
logarithme de $3/2$.  Il suffit d'accepter que l'idée d'additionner
$0.58496$ est intéressante, par paresse, ou comme un argument
d'autorité, ou en s'efforçant de comprendre les logarithmes, chacun
s'y retrouvera selon sa personnalité.  Ensuite, avec simplement
quelques figures, on trouve une explication aux gammes pentatoniques,
diatoniques et chromatiques, à la distinction entre les demi-tons
chromatiques et diatoniques, et à d'autres choses que je n'ai peut-être
pas vues, ce qui n'est déjà pas si mal.

Une nouvelle manière de comprendre les constructions du passé pourrait
être un point de départ de celles du futur.  J'ai à peine esquissé cet
aspect en proposant une construction basée sur l'intervalle de tierce
majeure pure.  Cette construction, et toutes celles qu'on pourrait
faire sur le même principe, ne nécessitent pas vraiment le théorème
des trois distances (encore une fois, c'est un cadre commode de
représentation) et les plus naturelles, ou disons les moins
artificielles, sont très probablement déjà connues, bien que je n'aie
pas trouvé de sources à leur sujet.  Poursuivre dans cette voie et en
aborder le côté réellement musical se heurtera nécessairement aux
habitudes culturelles: des mouvements conjoints dans des gammes
étranges sonneront bizarrement à nos oreilles, il y a peu de risque à
le parier.

\'Evidemment, il faudrait être bien naïf pour penser que l'application
de quelques formules simples, même issues de mathématiques récentes et
enrobées d'un vocabulaire abscons, produira un résultat artistique
intéressant, surtout après des millénaires de culture musicale de
plusieurs civilisations.  \c Ca n'en reste pas moins amusant, et la
naïveté est parfois une qualité chez le scientifique ou l'artiste,
alors pourquoi ne pas essayer?  Le risque est surtout de parvenir à la
conclusion que le système classique n'a pas survécu jusqu'à nous sans
de bonnes raisons et possède des qualités intrinsèques difficiles à
dépasser, ou même à imiter.  Pour illustrer ce type de circularités
vaines, je me permets pour conclure de citer la cinquième \oe uvre au
catalogue des écrits de Pierre Ménard, l'immortel auteur du Quichotte,
telle que relatée par Jorge Luis Borges:

{\it \og Un article technique sur la possibilité d'enrichir le jeu d'échecs
en éliminant un des pions de la tour. Ménard propose, recommande,
discute et finit par rejeter cette innovation\fg{}.}

\section*{Remerciements}

Je remercie Guillaume Hanrot pour m'avoir donné l'idée de ce travail,
et Claire Lapalu de l'ENM pour l'avoir validée.  Merci à Valérie
Berthé qui m'a donné de nombreuses pistes de travail m'ayant
finalement conduit à la thèse de Norman Carey. Je remercie ce dernier
pour m'avoir aimablement communiqué sa thèse.  Et merci à ma chère épouse
Christelle Petit d'avoir relu le document.


\end{document}